\documentclass[12pt]{article}
\newcommand{\bee}{\begin{eqnarray*}}
\newcommand{\ene}{\end{eqnarray*}}
\newcommand{\beeq}{\begin{equation}}
\newcommand{\eneq}{\end{equation}}
\newtheorem{lem}{Lemma}[section]
\newcommand{\bel}{\begin{lem}}
\newcommand{\enl}{\end{lem}}
\newtheorem{defi}{Definition}[section]
\newcommand{\bef}{\begin{defi}}
\newcommand{\enf}{\end{defi}}
\newtheorem{exap}{Example}[section]
\newcommand{\beex}{\begin{exap}}
\newcommand{\enex}{\end{exap}}
\newtheorem{theo}{Theorem}[section]
\newcommand{\beth}{\begin{theo}}
\newcommand{\enth}{\end{theo}}
\newtheorem{prop}{Proposition}[section]
\newcommand{\bep}{\begin{prop}}
\newcommand{\enp}{\end{prop}}
\newtheorem{cor}{Corollary}[section]
\newcommand{\bec}{\begin{cor}}
\newcommand{\enc}{\end{cor}}
\newtheorem{rem}{Remark}[section]
\newcommand{\ber}{\begin{rem}}
\newcommand{\enr}{\end{rem}}

\usepackage{amsmath}
\usepackage{graphicx}
\usepackage{latexsym}
\begin{document}
\title { 
STATISTICAL INFERENCE\\ WITH \\
DATA  AUGMENTATION\\
AND\\
PARAMETER EXPANSION
}
 \author {Yannis G. Yatracos\\
Faculty of Communication and Media Studies\\
Cyprus University of Technology}
\maketitle
\date{}

\maketitle
\date{}
\vspace{1.6in}

\pagebreak

\begin{center} \vspace{0.05in} {\large Summary} \end{center}
\begin{center}
\parbox{4.8in}
{\quad Statistical pragmatism 
embraces 
all 
efficient methods  in statistical inference. 
Augmentation of the collected data  is used herein  to
obtain  {\em  representative population  information}
from a {\em large} 
class of non-representative population's units.
Parameter expansion of a probability model  is shown
to reduce the 
upper bound on the sum of error probabilities  for a test of 
simple hypotheses, and a measure, $R,$ is proposed  for the effect  
of activating additional component(s) in the   
sufficient statistic.
\\

}
\end{center}

\vspace{.05 in}

\bigskip
{\it Some key words:} \quad Collected Data Augmentation;  Parameter Expansion; Representative Information;  Statistical Pragmatism

\mbox{ } 
\vspace{.05in}

\pagebreak

\section {Introduction}

\quad In recent years, response rates in polls have decreased further  and  a representative sample may not
 always 
become available. For example,
using random digit dialing (RRD)
the response rate
decreased from 36\% in 1997
to 9\% in 2012 (Kohut
{\em et al.} 2012). According to Wang {\em et al.} (2015),  it is convenient and cost-effective to collect
a very large non-representative sample via online surveys and obtain with 
statistical adjustments accurate  
election forecasts, on par with those based on traditional representative polls.

We study  initially  the problem of obtaining 
{\em  representative information}
 for the population,
from a large number of non-representative population's  units
with a common attribute, ${\cal A}.$  
 This attribute could be, e.g., 
an account in Facebook or   following  a celebrity in a Social Network;
the latter  occured in practice with Social Voting Advice Applications (Katakis {\em et al.}, 2014).

Under assumptions occuring also in practice, it is expected that units with  attribute  ${\cal A}$ 
can provide each additional, accurate  information for one of the remaining  units  in the population without attribute ${\cal A}.$ 
Due to the large number of units with attribute  ${\cal A},$  the so-augmented
data-information 
from all strata 
can be used to obtain
 information equivalent to that from  a representative sample (Kruskal and  Mosteller, 1979).

The second problem studied is  model parameter expansion (PX), which   is shown to
reduce the
 upper bound on the 
sum of error probabilities,
when testing two simple hypotheses with a  test introduced by Kraft (1955).
 The 
 proof  confirms that parameter expansion 
 is  ``activating'' a sufficient statistic with additional component(s)
(Rubin, 1997) and the effect of the activation  is clarified 
with the introduced  measure, $R,$ obtained with a one-parameter expansion.
 These results  explain
 why 
the PX-EM algorithm (Liu {\em et al.},1998)  converges faster than the 
EM algorithm (Dempster {\em et al.},
1977) and its many variations.


Fisher (1922) introduced  in statistical inference the use of  a 
model that is not updated.
Essentially all models are wrong but some are very useful (attributed to
George Box in Rubin, 2005). The need for improvement of estimation
procedures led statisticians to relax the use of the 
``one and only'' assumed model by adopting,
for example,
the Bayesian model averaging approach 
(see e.g.,  Hoeting {\it et al.}, 1999). 
Expansion of  a probability model or
the augmentation of collected data 
improve,  respectively, 
the
data's fit 
and 
 the  estimates  of  the  model's parameters.
Artificial data augmentation
has been used
in missing value problems 
(e.g. see 
Rubin, 1987), to improve the convergence of the EM algorithm (see, e.g., Meng and van Dyk, 1997)
and
to reduce the mean squared error of $U$-statistics,
in particular the unbiased estimates of
variance and covariance 
 (Yatracos, 2005). 
 Parameter expansion  (Meng and van Dyk, 1997, Liu {\em et al.}, 1998) and  model updated maximum likelihood estimate (MUMLE, Yatracos, 2015)
are also
examples of deviations from the  one-model approach. 
 
Kass (2011)
introduced statistical pragmatism. It is  a new philosophy
 which is inclusive of the Bayesian, frequentist and all other efficient approaches in inference.  
 Statistical pragmatism emphasizes the 
assumptions that connect statistical models with the observed data
and it  is implemented herein with  parameter expansion and   augmentation of the collected data.
The additional contribution of the  latter is that: {\em a)} it introduces a new component in the statistical inference set-up: each unit in the population provides, in addition, data 
for other units, and {\em b)} it includes  as goal to obtain  representative population information when  representative units/sample
are not available.

\section{Representative Information from Non-Representative Units with Augmentation of the Collected Data}


\quad Accessible units in the population with attribute ${\cal A}$  responding to a questionnaire are not representative, e.g., when their age is between 25 and 35,  and the collected information is not comparable to that of a random sample.
Assumptions are made for  attribute ${\cal A}$  and the population, allowing
to obtain  representative population  information from these non-representative units. 

{\em Assumption 1-Common Attribute:}  Attribute ${\cal A}$ is {\em common} in the population:
the number of respondents with attribute 
${\cal A}$ is large,  each respondent has   in its immediate environment one associated unit
without  attribute ${\cal A},$   and collectively the associated units come from all strata. 

{\em Assumption 2- Effective Information:}  
A  large percentage of 
units with attribute ${\cal A}$
have each  accurate, 
 questionnaire related
 information,   
for its associated  unit
{\em without} attribute ${\cal A}.$   

\bef
\label{eq:repinfo}
Representative population information
 is equivalent to 
that obtained from a representative sample from  the population.
\enf
\quad The Common Attribute and the Efficient Information  assumptions  guarantee that a large number of units  responding to the questionnaire  will provide each  information for its associated unit, thus information will be obtained  from all  strata. A  representative sample, e.g. random sample, from the so-obtained information will provide  representative information  for the population. 

The degree of accuracy of the  information given from  the respondent 
for a  unit  without attribute ${\cal A}$  may be clarified from the respondent's answers to the questionnaire.
Given the large number of respondents, only the most accurate information for  associated  units without attribute ${\cal A}$  will be used.
The additional information  provided by some of the units may introduce bias, thus bias reduction methods will guarantee the best result.

The proposed method of data augmentation will be used in   Voting Advice Applications before the general  elections in Spain (December 20, 2015). Note that
 collected data augmentation.differs from the notion of data augmentation introduced  in Tanner and Wong (1987) and  Gelman (2004) .




\section{Parameter Expansion in Statistical Procedures}

\quad For the EM-algorithm and its variations, the observed data model
$f(x_{obs}|\theta)$ and the augmented (called also complete) data
model $f(x_{com}|\theta)$ have the same parameter $\theta.$ For the 
PX-EM algorithm 
(Liu {\it et al.},1998), $f(x_{com}|\theta)$ 
contains another parameter with known value $\eta_0$
and is expanded to
a larger model $f_X(x_{com}|(\theta_{*}, \eta))$ with $\theta_{*}$
playing the role of $\theta.$ 
The complete data model is preserved when 
$\eta=\eta_0,$ that is,
$$f_X(x_{com}|(\theta_*,\eta_0))=f(x_{com}|\theta=\theta_*).$$
In several examples in the same paper it is observed that the PX-EM algorithm 
is faster than the EM-algorithm and its variations.

Testing simple hypotheses is an elementary but fundamental 
 problem in statistical inference.
For example, families of tests of simple hypotheses allow to obtain
a consistent estimate, with calculation of rates of convergence
of this estimate in Hellinger distance (LeCam, 1986).
Improved tests 
will 
increase the accuracy of the so-obtained estimate.
It will be shown 
that a test of simple hypotheses
is improved 
with parameter expansion.

In the sequel,  omitted domains of integration are determined by the integrands-densities.


\bef
 For densities $f, g$  defined on $R$ the
Hellinger distance, $H(f,g),$ is 
given by 
$$H^2(f,g) = \int
[f^{1/2}(x)-g^{1/2}(x)]^2 dx.$$
The affinity of $f, g$ is
\[
 \rho(f,g)= \int
 f^{1/2}(x)g^{1/2}(x) dx. 
\]
\enf
For $H(f,g)$
and $\rho$ 
 it holds: 
\begin{equation}
H^2(f,g) = 2(1- \rho(f,g)) \le 2
; 
\label{eq:hellprop}
\end{equation}
$H^2(f,g)=2 $ if and only if
$\rho(f,g)=0,$ 
i.e., if 
$f(x)g(x)=0 $ almost surely.

For a sample $X_1,\ldots,X_n$ with density either $f^n$ or $g^n,$
the larger $H^2(f^n,g^n)$ is
 the easier is
to determine either  the true density of the data
or, in parametric models, the 
true parameter.
Consistent testing
is guaranteed because  $\lim_{n \longrightarrow 
+\infty}H^2(f^n, g^n)=2.$

Assume that  sample $X_1,\ldots, X_n$   has density $f(x|\theta)$
and that a different  model parameter, $\eta,$ with known value $\eta_0$ is
already included in the model. 
It is shown that the upper  bound on the sum of error probabilities of a consistent test introduced by Kraft (1955) 
for  hypotheses
\begin{equation}
H_0: \theta=\theta_0 \mbox{  and  } H_1: \theta=\theta_1,
\label{ht1}
\end{equation}
is reduced with parameter expansion.

Let $t_{n,1}=t_1$ be the sufficient statistic for 
$\theta$ with density $g(t_1|\theta).$
Kraft (1955) provided the consistent test 
\begin{equation}
\phi_n=I(\frac{\sqrt{g(t_1|\theta_1)}}{\sqrt{g(t_1|\theta_0)}}>1);
\label{eq:kraft}
\end{equation}
$I$ denotes the indicator function. 
Note that the criticisms of the classical $\alpha$-level test for the preferential treatment of $H_0$  and the predetermined level of the test
do not hold for test (\ref{eq:kraft}).

The error probabilities when using  $\phi_n$ are
$$
E_{H_0}\phi_n \le  \int_{\sqrt{g(t_1|\theta_1)}>\sqrt{g(t_1|\theta_0)}} \sqrt{g(t_1|\theta_1)g(t_1|\theta_0)} dt_1,
$$
$$
E_{H_1}(1-\phi_n) \le \int_{\sqrt{g(t_1|\theta_0)}\ge \sqrt{g(t_1|\theta_1)}}  \sqrt{g(t_1|\theta_1)g(t_1|\theta_0)} dt_1.
$$
Then, the sum of error  probabilities
\begin{equation}
\label{eq:sumprob1}
E_{H_0}\phi_n+E_{H_1}(1-\phi_n) \le \int \sqrt{g(t_1|\theta_1)g(t_1|\theta_0)} dt_1.
\end{equation}

Lower values of the affinity $ \int \sqrt{g(t_1|\theta_1)g(t_1|\theta_0)} dt_1$
indicate an increase in the separation of the densities
$g(t_1|\theta_0)$ and $g(t_1|\theta_1),$ making easier to distinguish between the two hypotheses.

Consider the expanded  model $f(x|\theta, \eta)$
with the new sufficient statistics $(t_1, t_{2,n}=t_2)$ that have joint density
\begin{equation}
h(t_1, t_2| \theta, \eta)=g(t_1|\theta,\eta) \tilde g(t_2|t_1,\theta,\eta).
\label{eq:likextmod}
\end{equation}
Note that $g(t_1|\theta,\eta=\eta_0)=g(t_1|\theta)$ and that  $\tilde g$ is the conditional density of $t_2$ given $t_1.$

For the expanded model (\ref{eq:likextmod}) and the 
 hypotheses testing problem 
\begin{equation}
H'_0: \theta=\theta_0, \eta=\eta_0 \mbox{  and  } H'_1: \theta=\theta_1,
\eta=\eta_0, 
\label{ht2}
\end{equation}
statistic   $t_2$  is ``activated'' and
a consistent test for these hypotheses similar to $\phi_n$  is
$$\psi_n=I(\frac{\sqrt{g(t_1|\theta_1)\tilde g(t_2|t_1,\theta_1, \eta_0)}}
{\sqrt{g(t_1|\theta_0)\tilde g(t_2|t_1,\theta_0, \eta_0)}}>1).$$
For the sum of error probabilities 
of $\psi_n$  it holds
\begin{equation}
E_{H'_0}\psi_n + E_{H'_1}(1-\psi_n)
\le \int \int \sqrt{g(t_1|\theta_1, \eta_0)g(t_1|\theta_0, \eta_0)} 
\sqrt{\tilde g(t_2|t_1,\theta_1, \eta_0) 
\tilde g(t_2|t_1,\theta_0, \eta_0)} dt_2 dt_1. 
\label{type12e}
\end{equation}

When $\tilde g$ depends on $\theta$ and $t_1$ is fixed, it follows that  $\tilde g(t_2|t_1,\theta_1, \eta_0)$ 
and
$\tilde g(t_2|t_1,\theta_0, \eta_0)$ are not equal a.s. $t_2$ and from
Cauchy-Schwarz inequality for integrals
\begin{equation}
\int \sqrt{\tilde g(t_2|t_1,\theta_1, \eta_0)}
\sqrt{\tilde g(t_2|t_1,\theta_0, \eta_0)} dt_2<1,
\label{eq:Rubin}
\end{equation}
for every   $t_1.$

Using Fubini's theorem in the right side of  (\ref{type12e}) and (\ref{eq:Rubin}) it follows that
\begin{equation}
\int \int \sqrt{g(t_1|\theta_1, \eta_0)g(t_1|\theta_0, \eta_0)} 
\sqrt{\tilde g(t_2|t_1,\theta_1, \eta_0) 
\tilde g(t_2|t_1,\theta_0, \eta_0)} dt_2 dt_1
<
 \int \sqrt{g(t_1|\theta_1)g(t_1|\theta_0)} dt_1.
\label{type12e2} 
\end{equation}
 
From (\ref{eq:sumprob1}),(\ref{eq:Rubin})  and (\ref{type12e2})  it follows that the upper bound on the sum of error probabilities  in  (\ref{type12e}) for the expanded  model is smaller than that of the original model.

\bef
\label{def:Rubin}  The effect of  activating  component  $t_2$  in the sufficient statistic of the
expanded  model (\ref{eq:likextmod}) with respect to model (\ref{ht1})  is measured by the difference of the affinities,
\begin{equation}
\label{eq:Rubin1}
R=\int  \sqrt{g(t_1|\theta_1)g(t_1|\theta_0)} dt_1 -  \int \int  \sqrt{h(t_1, t_2|\theta_1, \eta_0)h(t_1,t_2|\theta_0, \eta_0)} dt_1 dt_2;
 \end{equation}
R is proportional to the difference of Hellinger distances of the  models in  (\ref{ht1}) and in   (\ref{eq:likextmod}).
 \enf





 
\begin{center}
{\bf Acknowledgments}
\end{center}

Many thanks are due to my colleague  Dr. Nicolas Tsapatsoulis and Dr. Fernando Mendez for introducing me to the theory
of Social Voting Advice Applications  and to  practical statistical  problems remaining to be solved.
Very many, special
thanks are due to the Department of Statistics and Applied Probability, National University of 
Singapore, for the  warm hospitality during my summer visits when most of the results were obtained. This research was partially supported by a grant from  Cyprus University of Technology.

\end{document}